\documentstyle[11pt,leqno]{article}
\input amssym.def
\catcode`\û=\active
\catcode`\…=\active
\catcode`\Ù=\active
\catcode`\–=\active
\catcode`\ð=\active
\catcode`\š=\active
\catcode`\ö=\active
\catcode`\û=\active
\letû=\ss
\def…{\"A}
\defÙ{\"a}
\def–{\"Oe}
\defð{\"o}
\defš{\"U}
\defö{\"u}
\title{Holomorphic spinors and the Dirac equation.}
\author{K.-D. Kirchberg, Berlin}
\date{\today}

\begin{document}

\maketitle

\begin{abstract} 
A closed spin KÙhler manifold of positive scalar curvature with smallest possible first eigenvalue of the Dirac operator is characterized by holomorphic spinors. It is shown that on any spin KÙhler-Einstein manifold each holomorphic spinor is a finite sum of eigenspinors of the square of the Dirac operator. Vanishing theorems for holomorphic spinors are proved.
\end{abstract}

\setcounter{section}{-1}
\section{Introduction}

Let $M^{2m}$ be a spin  KÙhler manifold of complex dimension $m$ and scalar curvature $R$. If $M^{2m}$ is closed and $R_0 := min (R) >0$, then it is known from [9] that the first eigenvalue $\lambda_1$ of the Dirac operator $D$ satisfies the inequality

\[ \lambda^2_1 \ge \left\{ \begin{array}{l} \displaystyle{\frac{m}{m-1} \cdot \frac{R_0}{4} \, \, , \quad m \, \, \mbox{even}}\\
\mbox{}\\
\displaystyle{\frac{m+1}{m} \cdot \frac{R_0}{4} \, \, , \quad m \, \, \mbox{odd} .} \end{array} \right. \]

$M^{2m}$ is called to be a limiting manifold iff this inequality is an equality. There are good descriptions of limiting manifolds for odd $m$ (see [8], [10], [11], [16]). The most important  result in this direction was proved by A. Moroianu which says that each limiting manifold of odd complex dimension $m \ge 3$ is the twistor space of a quaternionic KÙhler manifold of positive scalar curvature. By Theorem 6 in [11], this implies that limiting manifolds of complex dimension $m=4k+3$ are just the closed KÙhler-Einstein manifolds carrying a complex contact structure. Moreover, the complex projective spaces ${\Bbb C} P^m$ are the only limiting manifolds for $m=4k+1$. In this paper we add a holomorphic characterization of limiting manifolds: If $m$ is odd, then we prove that a closed spin KÙhler manifold $M^{2m}$ of positive scalar curvature is a limiting one iff $M^{2m}$ is Einstein and the bundle $\Lambda^{0, \frac{m+1}{2}} \otimes \sqrt{\Lambda^{m,0}}$ admits a holomorphic section, where $\Lambda^{m,0}$ is the canonical bundle and $\sqrt{\Lambda^{m,0}}$ the spin structure of $M^{2m}$. By Moroianu's theorem, this is an assertion on the twistor spaces of quaternionic KÙhler manifolds of positive scalar curvature. In case where $m$ is even we obtain an analogous description of limiting manifolds using a result of A. Moroianu concerning the eigenvalues of the Ricci tensor (see [16]). We recall that the situation is not clear for even $m \ge 4$. In 1990 A. Lichnerowicz announced a theorem (see [13]) which asserts that in this case each limiting  manifold is the  product of the flat torus $T^2$ by a limiting manifold of complex dimension $m-1$. But up to now there is no proof of this theorem. Thus, for even $m$, the only classification result is that of Th. Friedrich in case $m=2$ (see [2]). The holomorphic characterization of limiting manifolds given here is a generalization of the holomorphic description in case $m=2$ that was used by Th. Friedrich to obtain his classification.\\
Our paper contains some other results, too. For example, we show that on any spin KÙhler-Einstein manifold $M^{2m}$ each holomorphic spinor is the sum of at most $m+1$ eigenspinors of the square of the Dirac operator. Moreover, we obtain vanishing theorems for holomorphic spinors using basic Weitzenbðck formulas of spin KÙhler geometry.\\
The author would like to thank Th. Friedrich and A. Moroianu for useful hints and discussions.\\

\section{Basic equations satisfied by holomorphic spinors.}

Let $M^{2m}$ be a spin KÙhler manifold of complex dimension $m$ with complex structure $J$, KÙhler metric $g$ and spinor bundle $S$. Then the KÙhler form $\Omega$ defines a canonical splitting

\begin{equation}
S=S_0 \oplus S_1 \oplus \cdots \oplus S_m
\end{equation}

into holomorphic subbundles $S_r \cong \Lambda^{0,r} \otimes S_o \,\, (r=0,1, \ldots,m,)$ of rank $\left( \begin{array}{c} m\\r \end{array} \right)$, where $S_0 = \sqrt{\Lambda^{m,0}}$ is the square root of the canonical bundle representing the spin structure of $M^{2m}$ (see [6] or [3], Section 3.4.). Considering $\Omega$ as an endomorphism of $S$ the action of $\Omega$ on $S_r$ is just the multiplication by $i(m-2r)$. Let $ p(X) := \frac{1}{2} (X-i J X)$ and $\bar{p} (X):= \frac{1}{2} (X+i JX)$ for any real vector field $X$ and let $\psi \in \Gamma (S_r)$. Then we have $p(X) \cdot \psi \in \Gamma (S_{r+1})$ and $\bar{p}(X) \cdot \psi \in \Gamma (S_{r-1})$. Furthermore, let $(X_1, \ldots , X_{2m})$ be any local frame of vector fields and $(\xi^1, \ldots , \xi^{2m})$ the corresponding coframe. Using the notations $g_{kl}:= g(X_n, X_l), (g^{kl}):=(g_{kl})^{-1}$, $X^k := g^{kl} X_l$ (Einstein' s convention of summation) and $g(X) := g(X , \cdot)$ the holomorphic and antiholomorphic part of the covariant derivative $\nabla {\psi}$ of any spinor $\psi \in \Gamma (S)$ are locally defined by

\begin{eqnarray}
\begin{array}{c}
\nabla^{1,0} \psi := \xi^{k} \otimes \nabla_{p(X_k)} \psi = g (\bar{p}(X^k)) \otimes \nabla_{p(X_k)} \psi , \\
\mbox{}\\
\nabla^{0,1} \psi := \xi^{k} \otimes \nabla_{\bar{p}(X_k)} \psi = g (p(X^k)) \otimes \nabla_{\bar{p}(X_k)} \psi  , \end{array}
\end{eqnarray}

respectively. This yields the decomposition

\begin{equation}
\nabla \psi = \nabla^{1,0} \psi + \nabla^{0,1} \psi
\end{equation}

with $\nabla^{1,0} \psi \in \Gamma (\Lambda^{1,0} \otimes S)$ and $\nabla^{0,1} \psi \in \Gamma (\Lambda^{0,1} \otimes S)$. $\psi$ is called to be holomorphic (antiholomorphic) iff $\nabla^{0,1} \psi=0 \,\, (\nabla^{1,0} \psi =0)$.\\

The Dirac operator $D$ and its KÙhler twist $\tilde{D}$ are locally given by 

\begin{equation}
D = X^k \cdot \nabla_{X_k} \quad , \quad \tilde{D} = J (X^k) \cdot \nabla_{X_k} . 
\end{equation}

In the following the operators $D_+$ and $D_-$ appear defined by $D_{\pm} := \frac{1}{2} (D \mp i \tilde{D})$. Then it holds

\begin{eqnarray}
\begin{array}{l}
D_+ =p(X^k) \cdot \nabla_{X_k} = p(X^k) \cdot \nabla_{\bar{p} (X_k)} = X^k \cdot \nabla_{\bar{p} (X_k)} , \\
\mbox{}\\
D_- =\bar{p}(X^k) \cdot \nabla_{X_k} = \bar{p}(X^k) \cdot \nabla_{{p} (X_k)} = X^k \cdot \nabla_{{p} (X_k)}
\end{array}
\end{eqnarray}

and there are the operator equations

\begin{equation}
D= D_+ + D_- \, \, , 
\end{equation}

\begin{equation}
D_+^2 =D_-^2 =0 , 
\end{equation}

\begin{equation}
D^2 =D_+ D_- + D_- D_+ . 
\end{equation}

Moreover, for $\psi \in \Gamma (S_r)$, we have $D_{\pm} \psi \in \Gamma (S_{r \pm 1})$.\\

\underline{Proposition 1:} {\it Let $\psi \in \Gamma (S)$ be any holomorphic (antiholomorphic) spinor. Then $\psi$ satisfies the equation

\begin{equation}
D_+ \psi = 0 \quad \quad (D_- \psi =0)
\end{equation}

and, moreover, the equation

\begin{equation}
\begin{array}{c}
\nabla_{\bar{p}(X)} D_- \psi + \frac{1}{2} \bar{p}(Ric (X)) \cdot \psi =0\\
\mbox{}\\
(\nabla_{{p}(X)} D_+ \psi + \frac{1}{2} {p}(Ric (X)) \cdot \psi =0)\\
\end{array}
\end{equation}

for each real vector field $X$, where $Ric$ is the Ricci tensor.}\\

\underline{Proof:} By definition, $\psi \in \Gamma (S)$ is holomorphic (antiholomorphic) iff there is the equation

\begin{equation}
\nabla_{\bar{p} (X)} \psi =0 \quad \quad (\nabla_{p(X)} \psi=0) 
\end{equation}

for each real vector field $X$. Thus, for example, we have $0= X^k \cdot \nabla_{\bar{p}(X_k)} \psi =D_+ \psi$ and hence (9). For any $\varphi \in \Gamma (S)$ and any $X \in \Gamma (TM^{2m})$, there is the well-known relation

\begin{equation}
X^k \cdot C(X_k, X) \varphi = \frac{1}{2} Ric (X) \cdot \varphi , 
\end{equation}

where $C$ is the curvature tensor of $S$. Hence, since $M^{2m}$ is KÙhler, it holds

\begin{eqnarray*}
\frac{1}{2} \bar{p} (Ric(X)) \cdot \varphi &=& \frac{1}{2} Ric(\bar{p}(X)) \cdot \varphi = \frac{1}{2} X^k \cdot C(X_k, \bar{p}(X)) \varphi =\\
&=& X^k \cdot C(p(X_k), \bar{p}(X)) \varphi = \bar{p} (X^k) \cdot  C(p(X_k), \bar{p}(X)) \varphi =\\
&=& \bar{p} (X^k) \cdot C(X_k, \bar{p}(X)) \varphi .
\end{eqnarray*}

Thus, (12) implies in KÙhler case the relations

\begin{eqnarray}
\begin{array}{l}
\frac{1}{2} \bar{p} (Ric(X)) \cdot \varphi = \bar{p} (X^k) \cdot C(X_k, \bar{p} (X)) \varphi , \\
\mbox{}\\
\frac{1}{2} {p} (Ric(X)) \cdot \varphi = {p} (X^k) \cdot C(X_k, {p} (X)) \varphi . \\
\end{array}
\end{eqnarray}

Now, let $P \in M^{2m}$ be any point and let $(X_1, \ldots , X_{2m})$ be a frame on a neighbourhood of $P$ such that

\begin{equation}
(\nabla X_k )_P =0 \hspace{2cm} (k=1, \cdots, 2m) . 
\end{equation}

\medskip

Using (5), (11), (13) and (14) we have at $P$

\begin{eqnarray*}
\begin{array}{l}
\nabla_{\bar{p}(X)} D_- \psi + \frac{1}{2} \bar{p} (Ric(X)) \cdot \psi = \nabla_{\bar{p}(X)} D_- \psi + \bar{p} (X^k) \cdot  C (X_k, \bar{p}(X)) \psi =\\
\mbox{}\\
= \nabla_{\bar{p}(X)} (\bar{p}(X^k) \nabla_{X_k} \psi) + \bar{p} (X^k) (\nabla_{X_k} \nabla_{\bar{p}(X)} \psi - \nabla_{\bar{p}(X)} \nabla_{X_k} \psi - \nabla_{[X_k, \bar{p}(X)]} \psi ) =\\
\mbox{}\\
= - \bar{p} (X^k) \nabla_{[X_k, \bar{p}(X)]}\psi = - \bar{p} (X^k) \nabla_{\nabla_{X_k} \bar{p}(X)} \psi = - \bar{p} (X^k) \nabla_{\bar{p}(\nabla_{X_k}X)} \psi =0 .
\end{array}
\end{eqnarray*}

\medskip

This yields (10). \quad $\Box$

\bigskip

\underline{Corollary 2:} {\it Let $\psi \in \Gamma (S)$ be holomorphic (antiholomorphic). Then there is the equation

\begin{equation}
D_{\stackrel{+}{(-)}} D_{\stackrel{-}{(+)}} \, \, \psi = \frac{R}{4} \psi \begin{array}{c}+\\(-) \end{array} \frac{i}{2} \rho \cdot \psi , 
\end{equation}

where $R$ is the scalar curvature and $\rho$ the Ricci form.}\\

\medskip 

\underline{Proof:} Using (10) and the first one of the well-known relations

\begin{eqnarray}
\begin{array}{l}
X^k \cdot \bar{p} (Ric(X_k)) = - \frac{R}{2} - i \rho , \\
\mbox{}\\
X^k \cdot {p} (Ric(X_k)) = - \frac{R}{2} + i \rho , \\
\end{array}
\end{eqnarray}

we have

\begin{eqnarray*}
0 &=& X^k \cdot \nabla_{\bar{p}(X_k)} D_- \psi + \frac{1}{2} X^k \cdot \bar{p} (Ric(X_k)) \cdot \psi =\\
&=& D_+ D_- \psi - \frac{R}{4} \psi - \frac{i}{2} \rho \cdot \psi . \quad \quad \Box
\end{eqnarray*}

\medskip

Let $\nabla^* \nabla$ be the Bochner Laplacian on $\Gamma (S)$ locally given by

\begin{equation}
\nabla^* \nabla = - g^{kl} (\nabla_{X_k} \nabla_{X_l} - \nabla_{\nabla_{X_k} X_l} ) . 
\end{equation}

\medskip

Then there is the well-known relation (see [12])

\begin{equation}
\nabla^* \nabla = D^2 - \frac{R}{4}
\end{equation}

\medskip

Using (8), (9), (15) and (18) we immediately obtain\\

\underline{Corollary 3:} {\it If $\psi \in \Gamma (S)$ is holomorphic (antiholomorphic), then $\psi$ satisfies the equivalent equations}

\begin{equation}
D^2 \psi = \frac{R}{4} \psi \stackrel{+}{(-)} \frac{i}{2} \rho \cdot \psi ,
\end{equation}

\begin{equation}
\nabla^* \nabla \psi = \stackrel{+}{(-)} \frac{i}{2} \rho \cdot \psi .
\end{equation}

\section{Holomorphic spinors and limiting manifolds.}

\underline{Theorem 4:} {\it Let $M^{2m}$ be a closed spin KÙhler manifold of odd complex dimension $m$ with positive scalar curvature and spin structure $S_0= \sqrt{\Lambda^{m,0}}$. Then $M^{2m}$ is a limiting manifold iff $M^{2m}$ is Einstein and the bundle $\Lambda^{0, \frac{m+1}{2}} \otimes S_0$ admits a holomorphic section.}\\

\underline{Proof:} We recall that there are canonical unitary isomorphisms

\begin{equation}
\alpha_r : \Lambda^{0,r} \otimes S_0 \stackrel{\sim}{\longrightarrow} S_r \quad \quad \quad  (r=0,1, \ldots , m) 
\end{equation}

defined by $\alpha_r (\omega \otimes \psi_0):= 2^{- \frac{r}{2}} \omega \cdot \psi_0$ (see [7], Prop. 4). From Section 4 in [9] we know that limiting manifolds of odd complex dimension $m$ are Einstein. Moreover, there is a holomorphic eigenspinor $\psi \in \Gamma (S_{\frac{m+1}{2}})$ to the first eigenvalue $\lambda^2_1 = \frac{m+1}{4m} R $ of $D^2$. Thus, $\left( \alpha_{\frac{m+1}{2}} \right)^{-1} \circ \psi$ is a holomorphic section of $\Lambda^{0, \frac{m+1}{2}} \otimes S_0$. Conversely, let $M^{2m}$ be Einstein and let $\varphi \not\equiv 0$ be a holomorphic section of $\Lambda^{0, \frac{m+1}{2}} \otimes S_0$. Then $\psi := \alpha_{\frac{m+1}{2}} \circ \varphi \in \Gamma (S_{\frac{m+1}{2}})$ is holomorphic. Hence, using Corollary 3, the Einstein condition $\rho = \frac{R}{2m} \Omega$ and $\Omega \psi = i(m-2 \cdot \frac{m+1}{2}) \psi = - i \psi$ we have

\[ D^2 \psi= \frac{R}{4} \psi + \frac{i}{2} \rho \psi = \frac{R}{4} \psi + \frac{R}{4m} \psi = \frac{m+1}{4m} R \psi . \]

Thus, $M^{2m}$ is a limiting manifold. \hspace{1cm} $\Box$\\

\medskip

The corresponding holomorphic characterization of limiting manifolds in case of even complex dimension $m$ is more complicated since such manifolds are not Einstein for $m \ge 4$. In case $m=2$ the complete list of limiting manifolds is given by $S^2 \times S^2$ and $S^2 \times T^2$ (see [2]). In case $m=2l \ge 4$ it is known that the scalar curvature $R$ is a (positive) constant and that there is a holomorphic eigenspinor $\psi \in \Gamma (S_{\frac{m}{2}})$ to the first eigenvalue $\lambda^2_1 = \frac{m}{4(m-1)} R$ of $D^2$ which additionally satisfies the eigenvalue equation

\begin{equation}
\rho \cdot \psi = -i \frac{R}{2m-2} \psi
\end{equation}

(see [9], Section 4). A previous result is\\

\underline{Proposition 5:} \, \, \, {\it Let $M^{2m}$ be a closed spin KÙhler manifold of even complex dimension $m \ge 4$ with positive scalar curvature $R$. Then $M^{2m}$ is a limiting  manifold iff $R$ is constant and $S_{\frac{m}{2}} $ admits a holomorphic section satisfying equation (22).}\\

\underline{Proof:} Let  $0 \not\equiv \psi \in \Gamma (S_{\frac{m}{2}})$ be holomorphic and also a solution of equation (22) and let $R>0$ be constant. Then, by Corollary 3, it holds

\[ D^2 \psi = \frac{R}{4} \psi + \frac{i}{2} \rho \psi = \frac{R}{4} \psi + \frac{R}{4(m-1)} R \psi = \frac{m}{4(m-1)} R \psi .  \]

Hence, $M^{2m}$ is a limiting manifold. The converse is true by the preceeding remarks. \hspace{3cm} $\Box$\\

\medskip

In the following we replace the condition (22) by a more geometrical one. We use the result of A. Moroianu that the Ricci tensor of a limiting manifold $M^{2m}$ with $m=2l \ge 4$ has exactly the two constant eigenvalues $\frac{R}{2m-2}$ and $0$ of  multiplicity $2m-2$ and $2$, respectively (see [16]). This property is called the condition (Ric).\\

\underline{Lemma 6:} \, \, {\it Let $M^{2m}$ be a spin KÙhler manifold satisfying the condition (Ric) and let $r \in \{1, \ldots, m-1 \}$. Then there is an orthogonal splitting

\begin{equation}
S_r =S_{r,0} \oplus S_{r,1}
\end{equation}

into complex subbundles with $rank_{\Bbb C} (S_{r, \varepsilon})= \left( \begin{array}{c}m-1\\r - \varepsilon \end{array} \right)$ $(\varepsilon =0,1)$ such that $-i \rho$ acts on $S_{r, \varepsilon}$ as multiplication by $\frac{R}{2m-2} \cdot (m-1-2(r- \varepsilon))$.}\\

\underline{Proof:} Let  $P \in M^{2m}$ be any point. The Ricci form $\rho$ defines an endomorphism $\rho_P : S_P \to S_P$ being compatible with the splitting $S_P=(S_0)_P \oplus \cdots \oplus (S_m)_P$. \\

If we look at the spin representation, then there is an identification such that

\begin{equation}
\Omega_P = \sum\limits^{m}_{k=1} e_{2k-1} \cdot e_{2k} \quad , \quad \rho_P = \sum\limits^{m}_{k=1} \rho_k (P) e_{2k-1} \cdot e_{2k} , 
\end{equation}

where $\rho_1 (P), \ldots, \rho_m (P)$ are the eigenvalues of the Ricci tensor at $P$ and $(e_1, \cdots , e_{2m})$ is the canonical basis of ${\Bbb R}^{2m}$ corresponding to a suitable orthonormal frame of $T_P M^{2m}$. Let $\{u_{\varepsilon_1  \ldots \varepsilon_m} | \varepsilon_1, \cdots , \varepsilon_m \in \{1,-1\}\}$ be the standard basis of the spin module. Then there are the relations

\begin{equation}
e_{2k-1} \cdot e_{2k} \cdot u_{\varepsilon_1  \ldots \varepsilon_m} = i \varepsilon_k \cdot u_{\varepsilon_1  \ldots \varepsilon_m} \quad (k=1, \cdots , m) . 
\end{equation}

This yields

\begin{equation}
\rho_P \cdot u_{\varepsilon_1  \ldots \varepsilon_m} = i \left( \sum\limits^m_{k=1} \varepsilon_k \rho_k (P) \right) \cdot u_{\varepsilon_1  \ldots \varepsilon_m} . 
\end{equation}

In case of condition (Ric) we can assume that \\
$\rho_1 (P) = \cdots = \rho_{m-1}(P) = \frac{R}{2m-2}, \rho_m (P) =0$. Then (26) implies

\[ \rho_P \cdot u_{\varepsilon_1  \ldots \varepsilon_m} = i \frac{R}{2m-2} \left( \sum\limits^{m-1}_{k=1} \varepsilon_k \right) \cdot u_{\varepsilon_1  \ldots \varepsilon_m} . \]

Now, $(S_r)_P$ corresponds to the vector space over ${\Bbb C}$ spanned by all $u_{\varepsilon_1  \ldots \varepsilon_m}$ for which exactly $r$ of the $\varepsilon_k$ are equal to $-1$. Hence, corresponding to the two possibilities $\varepsilon_m =1$ and $\varepsilon_m =-1$ the endomorphisms $\rho_P$ restricted to $(S_r)_P$ has exactly the two eigenvalues $i \frac{R}{2m-2} (m-1-2r)$ and $i \frac{R}{2m-2} (m-1-2(r-1))$ of multiplicity $\left( \begin{array}{c}m-1\\r \end{array} \right)$ and $\left( \begin{array}{c} m-1\\r-1 \end{array} \right)$, respectively. \, \, \, $\Box$\\

\medskip

Lemma 6 shows that in case of condition (Ric) a spinor $\psi \in \Gamma (S_{\frac{m}{2}})$ satisfies the equation (22) iff $\psi$ is a section of the subbundle $S_{\frac{m}{2}, 0}$ of $S_{\frac{m}{2}}$. Now we will see how this subbundle can be constructed with help of the Ricci tensor. The condition (Ric) provides an orthogonal $J$-invariant decomposition

\begin{equation}
TM^{2m} = E \oplus F , 
\end{equation}

where the fibres of the subbundles $E$ and $F$ at any point $P \in M^{2m}$ are given by $E_P:= \ker (Ric_P - \frac{R}{2m-2})$ and $F_P:= \ker (Ric_P)$, respectively. This implies the decomposition

\begin{equation}
T^{0,1*} M^{2m} = E^{0,1*} \oplus F^{0,1*}
\end{equation}

with $rank_{\Bbb C} (E^{0,1*} )=m-1$ and $rank_{\Bbb C} (F^{0,1*})=1$. Thus, we have

\begin{eqnarray*}
\Lambda^{0,r} :&=& \Lambda^r (T^{0,1*} M^{2m}) = \Lambda^r (E^{0,1*} \oplus F^{0,1*}) =\\
&=& (\Lambda^r E^{0,1*}) \oplus ((\Lambda^{r-1} E^{0,1*}) \otimes F^{0,1*}) . 
\end{eqnarray*}

Hence, using the isomorphisms (21) we obtain

\begin{equation}
(\Lambda^r E^{0,1*}) \otimes S_0 \oplus (\Lambda^{r-1} E^{0,1*}) \otimes F^{0,1*} \otimes S_0 \cong S_r . 
\end{equation}

By construction, it holds

\begin{equation}
\begin{array}{c}
rank_{\Bbb C} (( \Lambda^r E^{0,1*}) \otimes S_0)= \left( \begin{array}{c} m-1\\r \end{array} \right) , 
\mbox{}\\
\mbox{}\\
rank_{\Bbb C} (( \Lambda^{r-1} E^{0,1*}) \otimes F^{0,1*} \otimes S_0)= \left( \begin{array}{c} m-1\\r-1 \end{array} \right) .
\end{array}
\end{equation}

\medskip

\underline{Lemma 7:} \, \, {\it Let $r \in \{1, \ldots, m-1 \}$. Then the isomorphism (29) induces isomorphisms}

\begin{equation}
\begin{array}{c}
(\Lambda^r E^{0,1*}) \otimes S_0 \cong S_{r,0} , \\
\mbox{}\\
(\Lambda^{r-1} E^{0,1*}) \otimes F^{0,1*} \otimes S_0 \cong S_{r,1} . 
\end{array}
\end{equation}

\medskip

\underline{Proof:} \, By Lemma 6 and (30), it is sufficient to show that $\alpha_r ((\Lambda^r E^{0,1*}) \otimes S_0)= S_{r,0}$. Hence, we have to prove that $\rho$ is the multiplication by $i \frac{R}{2m-2} (m-1-2r)$ on $\alpha_r ((\Lambda^r E^{0,1*}) \otimes S_0)$. For any local frame $(X_1, \ldots, X_n), n:=2m$, the Ricci form $\rho$ acts on $S$ by

\begin{equation}
\rho = \frac{1}{2} J (X^k) \cdot Ric (X_k)
\end{equation}

(see [9], (54)). Now, let $(X_1, \ldots, X_n)$ be orthonormal such that

\begin{equation}
J (X_{2k})=X_{2k-1} \quad \quad (k=1, \ldots, m) , 
\end{equation}

\begin{equation}
Ric (X_{2k})= \frac{R}{2m-2} X_{2k} \quad \quad (k=1, \ldots, m-1), \qquad Ric(X_n)=0 . 
\end{equation}

Then we have

\begin{equation}
\Omega = \frac{1}{2} \sum\limits^n_{k=1} J (X_k) \cdot X_k
\end{equation}

and, moreover, by (32) and (34),

\begin{equation}
\rho = \frac{R}{2(n-2)} \cdot \sum\limits^{n-2}_{k=1} J (X_k) \cdot X_k . 
\end{equation}

Let us consider the 2-form $\eta := \Omega - \frac{n-2}{R} \rho$. Then, by (33), (35) and (36), we obtain

\[ \eta = \frac{1}{2} (J (X_{n-1}) \cdot X_{n-1} + J (X_n) \cdot X_n) = \frac{1}{2} ( - X_n \cdot X_{n-1} + X_{n-1} \cdot X_n) = X_{n-1} \cdot X_n . \]

Using this we determine the action of $\eta_P$ on the space

\[ \alpha_r ((\Lambda^r E^{0,1*}) \otimes S_0 )_P = \{ \omega \cdot \psi_0 |\omega \in (\Lambda^r E^{0,1*})_P, \psi_0 \in (S_0)_P \} \]

for any $P \in M^{2m}$. We remark that $\omega \in (\Lambda^r E^{0,1*})_p$ implies $X_{n-1} \lrcorner \omega =X_n \lrcorner \omega=0$. Hence, using (25) and the properties of Clifford multiplication we calculate

\[ \eta \cdot (\omega \cdot \psi_0) =X_{n-1} \cdot X_n \cdot (\omega \cdot \psi_0)=(-1)^r X_{n-1} \cdot (\omega \cdot (X_n \cdot \psi_0))= \]

\[ = (-1)^r  \cdot (-1)^r \omega \cdot (X_{n-1} \cdot X_n \cdot \psi_0)= i \omega \cdot \psi_0 . \]

Thus, $\eta$ is the multiplication by $i$ on $\alpha_r ( \Lambda^r E^{0,1*} \otimes S_0) $. This proves the assertion since $\rho = \frac{R}{n-2} (\Omega - \eta)$ and $\Omega$ is the multiplication by $i (m-2r)$ on $S_r$. \hspace{1cm} $\Box$\\

\medskip

Let $M^{2m}$ be a spin KÙhler manifold of even complex dimension $m$ satisfying the condition (Ric). Then the subbundle

\[ (\Lambda^{\frac{m}{2}} E^{0,1*} ) \otimes S_0 \subset \Lambda^{0, \frac{m}{2}} \otimes S_0 \]

constructed with help of the Ricci tensor only is not holomorphic in general. (Clearly, it is holomorphic if Ric is parallel). But also in case of a limiting manifold we can not say up to now whether this bundle must be holomorphic or not.  By Proposition 5, Lemma 6 and Lemma 7 we immediately obtain \\

\underline{Theorem 8:} \, {\it Let $M^{2m}$ be a closed spin KÙhler manifold of even complex dimension $m \ge 4$ and positive scalar curvature. Then $M^{2m}$ is a limiting manifold iff $M^{2m}$ satisfies the condition (Ric) and $\Lambda^{0, \frac{m}{2}} \otimes S_0$ admits a holomorphic section which simultaneously is a section of the subbundle $(\Lambda^{\frac{m}{2}} E^{0,1*}) \otimes S_0$.}\\

\medskip

We remark that the conditon (Ric) holds obviously if $M^{2m} =N^{2m-2} \times T^2$, where $N^{2m-2}$ is a limiting manifold of odd complex dimension $m-1$ (see [9], Section 5). In these cases the Ricci tensor is always parallel. Hence, if the theorem of Lichnerowicz is valid, then the Ricci tensor must be parallel for each limiting manifold of even complex dimension $m \ge 4$. But Theorem 8 suggests the conjecture that there are examples with non-parallel Ricci tensor.\\

\section{Holomorphic spinors on Einstein manifolds.}

Let $j : S \to S$ be the $j$-structure of the spinor bundle $S$ which always exists in even real dimensions. We recall that $j$ is parallel and anti-linear, preserves the length of spinors, commutes with Clifford multiplication by real vectors and has the properties

\begin{equation} j S_r = S_{m-r} \hspace{2cm} (r= 0,1, \ldots, m) , 
\end{equation}

\begin{equation}
j^2 =(-1)^{\frac{m(m+1)}{2}} . 
\end{equation}

Then we have obviously

\begin{equation}
j \ker (\nabla^{1,0}) = \ker (\nabla^{0,1}) \quad , \quad  j \ker (\nabla^{0,1}) = \ker (\nabla^{1,0}) . 
\end{equation}

Furthermore, let  $\lambda$ be any eigenvalue of $D^2$ and let $E_{\lambda}(D^2)$ denote the corresponding eigenspace. Then the relation

\begin{equation}
j E_{\lambda} (D^2) = E_{\bar{\lambda}} (D^2)
\end{equation}

is also obvious. $E_{\lambda}^{1,0} (D^2):=E_{\lambda} (D^2) \cap \ker (\nabla^{0,1}) \left(E_{\lambda}^{0,1} (D^2):= E_{\lambda} (D^2) \cap \ker (\nabla^{1,0})\right)$ is  the corresponding holomorphic (antiholomorphic) eigenspace. Moreover,  we  say that $\lambda$ is holomorphic (antiholomorphic) iff $E^{1,0}_{\lambda} (D^2) \not= 0$ $(E^{0,1}_{\lambda} (D^2) \not= 0)$. From (39) and (40) we see that

\begin{equation}
j E^{1,0}_{\lambda} (D^2)= E^{0,1}_{\bar{\lambda}} (D^2) \quad , \quad j E^{0,1}_{\lambda} (D^2) = E^{1,0}_{\bar{\lambda}} (D^2) . 
\end{equation}

Thus, $\lambda$ is holomorphic iff $\bar{\lambda}$ is antiholomorphic.\\

\medskip

\underline{Proposition 9:} \, {\it Let $M^{2m}$  be any spin KÙhler-Einstein manifold. Then the sets of holomorphic eigenvalues of $D^2$ coincides with the set of antiholomorphic eigenvalues and is  contained in $\{ \frac{rR}{2m} \, \,  | \,  r=0,1,\ldots,m \}$. Moreover, it holds}

\begin{equation}
\begin{array}{l}
E^{1,0}_{\frac{rR}{2m}} (D^2) \subseteq \Gamma (S_r) , \\
\mbox{}\\
E^{0,1}_{\frac{rR}{2m}} (D^2) \subseteq \Gamma (S_{m-r}) .
\end{array} \hspace{2cm} (r=0,1, \ldots ,m)
\end{equation}

\medskip

\underline{Proof:} \, Let $\lambda$ be any holomorphic eigenvalue of $D^2$. Moreover, let $0 \not\equiv \psi \in E^{1,0}_{\lambda} (D^2)$ and $\psi = \psi_0 + \psi_1+ \cdots + \psi_m$ the decomposition according to (1). Using the equation (19) and the Einstein condition

\begin{equation}
\rho = \frac{R}{2m} \Omega
\end{equation}

we have

\[ 0= \lambda \psi - D^2 \psi = \lambda \psi - \frac{R}{4} \psi - \frac{i}{2} \rho \psi = \left( \lambda - \frac{R}{4}\right) \psi - \frac{R}{4m} i \Omega \psi = \]

\[ = \left( \lambda - \frac{R}{4} \right) \psi + \frac{R}{4m} \sum\limits^m_{s=0} (m-2s) \psi_s = \sum\limits^m_{s=0} \left(\lambda - \frac{sR}{2m}\right) \psi_s \]

and hence $(\lambda - \frac{sR}{2m}) \psi_s =0$ for $s=0,1, \ldots, m$. Since $\psi \not\equiv 0$ there is an $r$ with $\psi_r \not\equiv 0$. \\

This implies $\lambda = \frac{rR}{2m}$ and  $\psi_s =0$ for $s \not= r$. \hspace{1cm} $\Box$\\

\medskip

\underline{Proposition 10:} \, {\it  Let $M^{2m}$ be a spin KÙhler-Einstein manifold, let  $r \in \{0,1, \ldots,m \}$ and let $\psi \in \Gamma (S_r)$ be holomorphic (antiholomorphic). Then $\psi$ satisfies the eigenvalue equation}

\begin{equation}
D^2 \psi = \frac{rR}{2m} \psi \quad \quad (D^2 \psi = \frac{(m-r)R}{2m} \psi ) . 
\end{equation}

\medskip

\underline{Proof:} \, For example, let $\psi \in \Gamma (S_r)$ be holomorphic. Then the equations (19), (43) and $\Omega \psi = i (m-2r) \psi$ imply (44). \hspace{1cm} $\Box$.\\

\medskip

\underline{Theorem 11:} \, {\it If $M^{2m}$ is a spin KÙhler-Einstein manifold of scalar curvature $R$, then there are the decompositions}

\begin{equation}
\begin{array}{l}
\ker (\nabla^{0,1} ) = \bigoplus\limits^m_{r=0} E^{1,0}_{\frac{rR}{2m}} (D^2) , \\
\mbox{}\\
\ker (\nabla^{1,0} ) = \bigoplus\limits^m_{r=0} E^{0,1}_{\frac{rR}{2m}} (D^2) . 
\end{array}
\end{equation}

\medskip

\underline{Proof:} \, Let $\psi \in \ker (\nabla^{0,1})$ and $\psi = \psi_0 + \psi_1 + \cdots + \psi_m$ the decomposition according to (1). Since the splitting (1) is parallel, $\psi$ is holomorphic iff each of its components $\psi_r$ is holomorphic. Thus, Proposition 10 implies $\psi_r \in E^{1,0}_{\frac{rR}{2m}} (D^2)$. \hspace{2cm} $\Box$\\

\medskip

From Proposition 10 and the structure (45) of the spaces of holomorphic or antiholomorphic spinors we immediately obtain\\

\underline{Theorem 12:} \, {\it Let $M^{2m}$ be a spin KÙhler-Einstein manifold and $R$  its scalar curvature. Then the following  holds: 
\begin{itemize}
\item[(i)] If $\frac{rR}{2m}$ is not an eigenvalue of $D^2$ for an $r \in \{0,1, \ldots,m\}$, then there is no holomorphic section in the bundle $S_r$  and no antiholomorphic section in $S_{m-r}$.
\item[(ii)] If $\frac{rR}{2m}$ is not an eigenvalue of $D^2$ for each $r \in \{ 0,1,\ldots,m\}$, then there are no holomorphic and no antiholomorphic spinors on $M^{2m}$.
\end{itemize}}

\section{Holomorphic spinors on closed manifolds.}

Let $Z,W$ be any complex vector  fields on a spin KÙhler manifold $M^{2m}$. Then we use the notation

\[ \nabla_{Z,W} := \nabla_Z \nabla_W - \nabla_{\nabla_{Z} W} \]

for the corresponding second order derivative of spinors. We consider the KÙhler-Bochner Laplacians on $\Gamma (S)$ locally defined  by

\begin{equation}
\begin{array}{l}
\nabla^{1,0*} \nabla^{1,0} := - g^{kl} \nabla_{\bar{p} (X_k) , p(X_l)} ,\\
\mbox{}\\
\nabla^{0,1*} \nabla^{0,1} := - g^{kl} \nabla_{p (X_k), \bar{p} (X_l)} .
\end{array}
\end{equation}

By definition, we have the inclusions

\begin{equation}
\ker (\nabla^{1,0} ) \subseteq \ker (\nabla^{1,0*} \nabla^{1,0}) \,\, , \, \, \ker (\nabla^{0,1}) \subseteq \ker (\nabla^{0,1*} \nabla^{0,1}) .
\end{equation}

\medskip

\underline{Proposition 13:} \, {\it There are the operator equations}

\begin{equation}
\begin{array}{l}
2 \nabla^{1,0*} \nabla^{1,0} = D^2 - \frac{R}{4} + \frac{i}{2} \rho ,\\
\mbox{}\\
2 \nabla^{0,1*} \nabla^{0,1} = D^2 - \frac{R}{4} - \frac{i}{2} \rho.
\end{array}
\end{equation}

\underline{Proof:} For example, we prove the first one of these equations. Let $P \in M^{2m}$ be any point and $(X_1, \ldots , X_{2m})$ a frame in a neighbourhood of $P$ with property $(\nabla X_k)_P=0$ for $k=1,\ldots,2m$. Using the formulas (5), (8), (13) and (16) we have at the point $P$

\begin{eqnarray*}
2 \nabla^{1,0*} \nabla^{1,0} &=& - 2 g^{kl} \nabla_{\bar{p} (X_k)} \nabla_{p (X_l)} = (X^k X^l + X^l X^k) \nabla_{\bar{p} (X_k)} \nabla_{p(X_l)} =\\
&=& D_+ D_- + X^l X^k (\nabla_{p(X_l)} \nabla_{\bar{p}(X_k)} - C(p(X_l), \bar{p} (X_k))) =\\
&=& D_+ D_- + D_- D_+ - X^l X^k C(p(X_l), \bar{p}(X_k))=\\
&=& D^2 + X^l p(X^k)C(X_k, p(X_l))=D^2 + \frac{1}{2} X^l p(Ric(X_l))=\\
&=& D^2 - \frac{R}{4} + \frac{i}{2} \rho .
\end{eqnarray*}
\mbox{} \hfill $\Box$\\

In 1979 M.L. Michelsohn proved Weitzenbðck formulas which are very similar to (48). M.L. Michelsohn also showed that $\nabla^{0,1*} \nabla^{0,1}$ and $\nabla^{1,0*} \nabla^{1,0}$ are non-negative, elliptic and formally self-adjoint differential operators (see [14], Prop. 7.2, 7.6). We prefer the equations (48) since the Ricci form $\rho$ enters here explicitely being more convenient for applications considered here. For example, we remark that Corollary 3 also follows from (47) and (48) immediately.\\

\underline{Theorem14:} \, {\it Let $M^{2m}$ be a closed spin KÙhler manifold. Then a spinor $\psi$ is holomorphic (antiholomorphic) iff $\psi$ satisfies equation (19) or (20).}\\

\underline{Proof:} \, In closed case the inclusions (47) are equalities (see [14], Prop. 7.2). $\Box$\\

\underline{Theorem 15:} \, {\it Let $M^{2m}$ be a closed spin KÙhler-Einstein manifold of scalar curvature $R$. Then we have the following:
\begin{itemize}
\item[(i)] For any $r \in \{0,1,\ldots,m\}$, $\psi \in \Gamma (S_r)$ is holomorphic (antiholomorphic) iff $\psi$ satisfies the eigenvalue equation (44).
\item[(ii)] In case $R>0$ the bundles $S_r$ with $r \le m/2 \,\, (r \ge m/2)$ do not admit any holomorphic (antiholomorphic) section. In case $R=0$ each holomorphic or antiholomorphic spinor is parallel.
\item[(iii)] It holds $\ker (\nabla^{0,1})= \bigoplus\limits_{r>m/2} E^{1,0}_{\frac{rR}{2m}} (D^2)$ and  $\ker (\nabla^{1,0})= \bigoplus\limits_{r>m/2} E^{0,1}_{\frac{rR}{2m}} (D^2)$ .
\end{itemize}} 

\underline{Proof:} In Einstein case the equations (19) and (44) are equaivalent  for $\psi \in \Gamma (S_r)$. Thus, assertion (i) follows from Theorem 14 immediately. Moreover, using the equation (20) and the Einstein condition (43) we see that $\psi \in \Gamma (S_r)$ is holomorphic (antiholomorphic) iff $\psi$ satisfies the equation

\[ \nabla^* \nabla \psi = \stackrel{-}{(+)} \frac{m-2r}{4m} R\psi \]

which implies

\[ \langle  \nabla^* \nabla \psi , \psi \rangle = \stackrel{-}{(+)} \frac{m-2r}{4m} R |\psi|^2 . \]

Integrating this equation we find

\[ || \nabla \psi ||^2 = \stackrel{-}{(+)} \frac{m-2r}{4m} R || \psi ||^2 . \]

Let $R>0$  and $\psi \not\equiv 0$. Then we obtain a contradiction for $m > 2r  \, \, (m<2r)$. The case $m=2r$ provides $\nabla \psi=0$. But it is known that the existence of a parallel spinor implies $Ric=0$ being a contradiction to our assumtion $R>0$. Finally, in case $R=0$ we always obtain $\nabla \psi=0$. This proves (ii). The assertion (iii) follows from (ii) and Theorem 11. \hspace{1cm} $\Box$\\

An essential generalization of Theorem 15, (ii) is \\

\underline{Theorem 16:} \, {\it Let $M^{2m}$ be a closed non-Ricci-flat spin KÙhler manifold and let $\rho_1 (P) \ge \rho_2 (P) \ge \cdots \ge \rho_m (P)$ denote the eigenvalues of Ric at $P \in M^{2m}$. Then the bundle $S_r$ $(r \in \{0,1,\ldots,m\})$ does not admit any holomorphic (antiholomorphic) section if  at each point $P$ the condition 

\begin{equation}
\rho_1 (P) + \cdots + \rho_r (P) \le \frac{1}{4} R(P) \quad \quad (\rho_{m-r+1} (P) + \cdots + \rho_m (P) \ge \frac{1}{4} R(P)) 
\end{equation}

is satisfied.}\\

\medskip

\underline{Proof:} Since $\rho_1 (P) + \cdots + \rho_m (P) = \frac{1}{2} R(P)$, we see from (26) that the set  of all eigenvalues of the endomorphism $-i \rho_P$ restricted to $(S_r)_P$ is given by 

\begin{equation}
 \left\{ \frac{1}{2} R(P) - 2 (\rho_{i_1} (P) + \cdots +  \rho_{i_r} (P)) | 1 \le i_1 < \cdots < i_r \le m \right\}.
\end{equation}

Thus, for any $\psi \in \Gamma (S_r)$, the condition (49) yields $\langle -i \rho \psi, \psi \rangle \ge 0$ $(\langle - i \rho \psi, \psi \rangle \le 0)$. On the other hand, if $\psi \in \Gamma (S_r)$ is holomorphic (antiholomorphic), then (20) implies the equation

\[ || \nabla \psi ||^2 = \stackrel{-}{(+)} \frac{1}{2} \int_{M^{2m}} \langle -i \rho \cdot \psi , \psi \rangle \]

which provides a contradiction for $\langle - i \rho \cdot \psi, \psi \rangle \ge 0$ \,\, $( \langle - i \rho \cdot \psi, \psi \rangle \le 0)$. \hspace{2cm} $\Box$\\

\medskip

We remark that the assertion of Theorem 16 is also valid if the condition (49) is not satisfied on a subset of $M^{2m}$ of measure zero. Moreover, in case $r=0$ the condition (49) simply reduces to $R \ge 0 \, \, (R \le 0)$. Hence, if $M^{2m}$ is closed and non-Ricci-flat, the bundle $S_0 \, \, (S_m)$ does not admit any holomorphic (antiholomorphic) section in case $R \ge 0 \, \, (R \le 0)$. Theorem 16 immediately yields\\

\underline{Corollary 17:} \, {\it Let $M^{2m}$ be a closed non-Ricci-flat spin KÙhler manifold of complex dimension $m \ge 3$ with scalar curvature $R \ge 0$ such that at each point $P \in M^{2m}$ the Ricci tensor has the eigenvalues $\frac{1}{2m-2} R(P)$ and 0 of multiplicity $2m-2$ and 0, respectively. Then the bundles $S_r$ with $r \le \frac{m-1}{2} \, \, ( r \ge \frac{m+1}{2} )$ do not admit any holomorphic (antiholomorphic) section.}\\

By Theorem 8 and Corollary 17, we obtain \\

\underline{Proposition 18:} \, {\it If $M^{2m}$ is a limiting manifold of even complex dimension $m \ge 4$, then the bundles $S_r$ with $r \le \frac{m-2}{2} \, \, (r \ge \frac{m+2}{2})$ do not admit any holomorphic (antiholomorphic) section.}\\

\end{document}